\documentclass[12pt,reqno]{amsart}
\usepackage{amsmath,amsfonts,amssymb,amsthm,amsxtra,amscd}
\begin{document}

\title[Stability of syzygy bundles]{On the stability of syzygy bundles}

\author[I. Coand\u{a}]{Iustin Coand\u{a}}
\address{Institute of Mathematics of the Romanian Academy, 
         P. O. Box 1-764, RO-014700, Bucharest, Romania}
\email{Iustin.Coanda@imar.ro}



\thanks{Partially supported by CNCSIS grant ID-PCE no.51/28.09.2007}

\begin{abstract}
We are concerned with the problem of the stability of the  
syzygy bundles associated to base point free 
vector spaces of forms of the same degree d on the projective space of 
dimension n. We deduce directly, from Mark Green's vanishing theorem for 
Koszul cohomology, that any such bundle is stable if his rank is sufficiently 
high. With a similar argument, we prove the semistability of a certain 
syzygy bundle on a general complete intersection of hypersurfaces of 
degree d in the projective space. This answers a question of H. Flenner 
(1984). We then give an elementary proof of H. Brenner's criterion of 
stability for monomial syzygy bundles, avoiding the use of Klyachko's 
results on toric vector bundles. We finally 
prove the existence of stable syzygy bundles defined by monomials of the 
same degree d, of any possible rank, for n at least 3. This extends the 
similar result proved, for n=2, by L. Costa, P. Macias Marques and 
R.M. Miro-Roig (2009). The extension to the case n at least 3 has been also, 
independently, obtained by P. Macias Marques in his thesis (2009).  
\end{abstract}

\maketitle

\section*{Introduction}

Let ${\mathbb P}^n$, $n \geq 1$, be the projective $n$-space over an 
algebraically closed field $k$ of {\it arbitrary characteristic}, let 
$S = k[X_0,\ldots ,X_n]$ be the homogeneous coordinate ring of 
${\mathbb P}^n$ and let $d \geq 1$ be an integer. A $k$-vector subspace 
$V$ of $S_d$ is called {\it base point free} (b.p.f., for short) if 
$\forall \, x \in {\mathbb P}^n$, $\exists \, f \in V$ such that 
$f(x) \neq 0$. In this case, we denote by $M_{d,V}$ the kernel of the 
evaluation epimorphism ${\mathcal O}_{{\mathbb P}^n}\otimes_kV \rightarrow 
{\mathcal O}_{{\mathbb P}^n}(d)$. We say that $V$ is {\it monomial} if it is 
generated by monomials in $S_d$. In this case, $V$ is b.p.f. if and only if 
it contains $X_0^d,\ldots ,X_n^d$. 

Consider, now, the polynomials $P_n,\, Q_{n-1} \in {\mathbb Q}[T]$ defined 
by: 
\[
P_n(T) = \frac{(T+1)\cdot \ldots \cdot (T+n)}{n!},\  
Q_{n-1}(T) = \frac{P_n(T) - 1}{T}\  . 
\]
As it is well known, $\text{dim}_k\, S_d = P_n(d)$, $\forall \, d \geq 0$. 

In this paper we prove several results concerning the stability of 
$M_{d,V}$. The first one is the following: 

\vskip3mm 
{\bf 1. Theorem.}\quad 
\textit{Let} $n \geq 2$, $d \geq 2$ \textit{and} $m$ \textit{be integers and 
let} $V$ \textit{be a b.p.f.}, $m$-\textit{dimensional subspace of} 
$S_d$. \textit{If} $P_n(d-1) + Q_{n-1}(d-1) < m \leq P_n(d)$ \textit{then} 
$M_{d,V}$ \textit{is stable and if} $m = P_n(d-1) + Q_{n-1}(d-1)$ 
\textit{then} $M_{d,V}$ \textit{is semistable}. 
\vskip3mm 

The notion of stability we use is that of {\it slope stability}. Its 
definition is recalled at the beginning of Section 1. If $m = P_n(d)$, i.e., 
if $V = S_d$, the semistability of $M_d := M_{d,V}$ was proved, in 
characteristic 0, by Flenner \cite{fle} and its stability by R. Paoletti 
\cite{pao}. In arbitrary characteristic, the semistability of $M_d$ was 
proved by Brenner \cite{bre}, Corollary 7.1 and, recently, V.B. Mehta gave 
another proof of this result in the appendix to the paper of Langer 
\cite{lan}. 

Our approach to the proof of Theorem 1 is quite different. We shall actually 
show that Theorem 1 is a direct consequence of Mark Green's vanishing 
theorem for Koszul cohomology \cite{gre}, Theorem 3.a.1. With the same method, 
we shall also prove the following result which answers a question of Flenner 
\cite{fle}, Remark 2.8., and leads to an improvement of the estimate in his 
restriction theorem: 

\vskip3mm 
{\bf 2. Proposition.}\quad 
\textit{Let} $n\geq 2$, $d \geq 2$ \textit{and} $1 \leq c \leq n-1$ 
\textit{be integers}, \textit{let} $Y \subset {\mathbb P}^n$ \textit{be a 
general complete intersection of} $c$ \textit{hypersurfaces of degree} 
$d$ \textit{and let} $N_{d,c}$ \textit{be the kernel of the evaluation 
epimorphism} ${\mathcal O}_Y\otimes_k\text{H}^0{\mathcal O}_Y(d) \rightarrow 
{\mathcal O}_Y(d)$. \textit{Then} $N_{d,c}$ \textit{is semistable with 
respect to} ${\mathcal O}_Y(1)$. 
\vskip3mm 

Next, we concentrate on the case where $V$ is monomial. We firstly give, 
in Section 2, an elementary, characteristic free proof of Brenner's criterion 
\cite{bre}, Theorem 6.3., of stability for monomial syzygy bundles. Brenner 
uses in his proof the results of Klyachko \cite{kly} on toric bundles. 
Klyachko developed his theory in characteristic 0, but M. Perling 
\cite{per} remarked that Klyachko's results are valid in arbitrary 
characteristic. We replace, in our proof, Klyachko's results by the fact that 
the Koszul complex defined by a set of monomials is a complex of 
${\mathbb N}^{n+1}$-graded $S$-modules.  

Using Brenner's criterion, Costa, Macias Marques and Mir\'{o}-Roig proved 
recently, in \cite{cmm}, in response to Question 7.8. from Brenner \cite{bre},  
the following: 

\vskip3mm 
{\bf 3. Theorem.}\quad (Costa, Macias Marques, Mir\'{o}-Roig, \cite{cmm})\\ 
\hspace*{3mm} 
\textit{Assume} $n = 2$. \textit{Then, for every integers} $d \geq 1$ 
\textit{and} $3 \leq m \leq P_2(d)$, \textit{there exists an}   
$m$-\textit{dimensional}, \textit{b.p.f. monomial subspace} 
$V$ \textit{of} $S_d$ \textit{such that} $M_{d,V}$ \textit{is stable}, 
\textit{except for the case where} $d = 2$ \textit{and} 
$m = 5$, \textit{when} $M_{d,V}$ \textit{is only semistable}. 
\vskip3mm 

We extend, in Section 3, this result to the case $n \geq 3$ in the following: 

\vskip3mm 
{\bf 4. Theorem.}\quad 
\textit{Let} $n \geq 3$, $d \geq 1$, \textit{and} $n+1 \leq m \leq P_n(d)$ 
\textit{be integers}. \textit{Then there exists an}  
$m$-\textit{dimensional}, \textit{b.p.f. monomial subspace} 
$V$ \textit{of} $S_d$ \textit{such that} $M_{d,V}$ \textit{is stable}. 
\vskip3mm 

We prove this result by induction on $n$ and on the degree $d$.  
The initial case $n = 2$ is, 
of course, the theorem of Costa, Macias Marques and Mir\'{o}-Roig.  
We show, however, that one can make 
the proof of Theorem 4 independent of the constructions from Costa et al. 
\cite{cmm}.  

As a consequence of the Theorems 3 and 4, it follows that if $n \geq 2$, 
$d \geq 1$ and $n+1 \leq m \leq P_n(d)$ are integers and if $V$ is a 
{\it general} b.p.f. subspace of $S_d$ then $M_{d,V}$ is stable, except for  
$n = 2$, $d = 2$ and $m = 5$. We study this exceptional case in Example 1.3. 
and show that, in this case too, $M_{d,V}$ is stable for a general $V$. 

After the completion of this paper, the author was 
informed by R.M. Mir\'{o}-Roig that the result stated in Theorem 4 was also,  
independently, obtained by P. Macias Marques in his thesis \cite{mma} using 
a different, more combinatorial, approach.  

\section{Applications of Green's vanishing theorem} 

Let $(X,{\mathcal O}_X(1))$ be an $n$-dimensional polarized smooth projective 
variety. If ${\mathcal F}$ is a torsion free coherent sheaf on $X$, the 
{\it degree} of ${\mathcal F}$ with respect to ${\mathcal O}_X(1)$ is 
$\text{deg}\, {\mathcal F} := (\text{det}\, {\mathcal F} \cdot 
{\mathcal O}_X(1)^{n-1})$ and its {\it slope} is $\mu ({\mathcal F}) := 
\text{deg}\, {\mathcal F}/\text{rk}\, {\mathcal F}$. ${\mathcal F}$ is 
called ({\it semi}){\it stable} with respect to ${\mathcal O}_X(1)$ if, 
for every coherent subsheaf ${\mathcal F}^{\prime}$ of ${\mathcal F}$ with 
$0 < \text{rk}\, {\mathcal F}^{\prime} < \text{rk}\, {\mathcal F}$, one has 
$\mu ({\mathcal F}^{\prime})\  (\leq ) < \mu ({\mathcal F})$. We recall 
that if ${\mathcal F}^{\prime \prime}$ is the kernel of the epimorphism 
${\mathcal F} \rightarrow ({\mathcal F}/{\mathcal F}^{\prime})/
({\mathcal F}/{\mathcal F}^{\prime})_{\text{tors}}$ then $\text{rk}\, 
{\mathcal F}^{\prime \prime} = \text{rk}\, {\mathcal F}^{\prime}$ and 
$\text{deg}\, {\mathcal F}^{\prime \prime} \geq  
\text{deg}\, {\mathcal F}^{\prime}$ hence, when verifying the (semi)stability 
condition, one may assume that, moreover, ${\mathcal F}/
{\mathcal F}^{\prime}$ is torsion free. We shall use the following obvious: 

\vskip3mm 
{\bf 1.1. Criterion of (semi)stability.} \quad 
\textit{Let} $(X,{\mathcal O}_X(1))$ \textit{be as above and let} $E$ 
\textit{be a vector bundle} (=\textit{locally free sheaf}) \textit{on} $X$. 
\textit{If}, \textit{for every} $r$ \textit{with} $0 < r < 
\text{rk}\, E$ \textit{and for every line bundle} $L$ 
\textit{on} $X$ \textit{with} $\mu ((\overset{r}{\wedge} E)\otimes L)\   
(<) \leq 0$, \textit{one has} $\text{H}^0((\overset{r}{\wedge} E)\otimes L) 
= 0$ \textit{then} $E$ \textit{is} (\textit{semi})\textit{stable with 
respect to} ${\mathcal O}_X(1)$. 
\vskip3mm 

For $(X,{\mathcal O}_X(1)) = ({\mathbb P}^n,{\mathcal O}_{{\mathbb P}^n}(1))$ 
the converse of the Criterion of semistability is also true, at least 
in characteristic 0. However, the converse of the Criterion of stability is 
not true. We shall see a counterexample in Example 1.3. below. 

We recall now, in a form which is more convenient for our purposes, 
Mark Green's 
vanishing theorem for Koszul cohomology. We also reproduce, for the reader's 
convenience, Green's elementary but ingenious argument. 

\vskip3mm 
{\bf 1.2. Lemma.} (Green \cite{gre}, 3.a.1.)\quad 
\textit{Let} $L$, $L^{\prime}$ \textit{be line bundles on a projective 
variety} $X$, $V$ \textit{a b.p.f. subspace of} $\text{H}^0(L)$ \textit{and} 
$M_{L,V}$ \textit{the kernel of the evaluation epimorphism} 
${\mathcal O}_X\otimes_kV \rightarrow L$.  
\textit{Then} $\text{H}^0((\overset{r}{\wedge} M_{L,V})\otimes 
L^{\prime}) = 0$ \textit{for} $r \geq \text{h}^0(L^{\prime}) := 
\text{dim}_k\, \text{H}^0(L^{\prime})$. 
\vskip3mm 

\begin{proof} 
Let us denote $M_{L,V}$ by $M$, $\text{dim}_k\, V$ by $m$ and 
$\text{h}^0(L^{\prime})$ by $l$. By hypothesis, $r \geq l$. 
We may assume, moreover,  
that $r \leq \text{rk}\, M = m-1$. Tensorizing by $L^{\prime}$ 
the exact sequence: 
\[
0\longrightarrow \overset{r}{\wedge} M\longrightarrow {\mathcal O}_X
\otimes_k\overset{r}{\wedge} V\overset{d}{\longrightarrow} L\otimes_k 
\overset{r-1}{\wedge} V        
\]
one deduces an exact sequence: 
\[
0\longrightarrow \text{H}^0(\overset{r}{\wedge} M\otimes L^{\prime}) 
\longrightarrow \text{H}^0(L^{\prime})\otimes \overset{r}{\wedge} V 
\overset{d}{\longrightarrow} \text{H}^0(L\otimes L^{\prime})\otimes 
\overset{r-1}{\wedge} V
\]
where $d$ is the Koszul differential. If  
$V^{\prime}$ be a non-zero subspace of $\text{H}^0(L^{\prime})$ and if 
$x \in X$ does not belong to the base locus of $V^{\prime}$ then: 
\[
\text{dim}_k\, \{g \in V^{\prime}\, \vert \, g(x) = 0\} = 
\text{dim}_k\, V^{\prime} - 1. 
\]
One deduces that if $x_1, \ldots ,x_m$ are general points of $X$ then: 

(i) \textit{If} $f \in V$ \textit{vanishes at} $x_1, \ldots ,x_m$ 
\textit{then} $f = 0$,\\ 
\hspace*{3mm} (ii) $\forall \, 1\leq i_1 < \cdots < i_l \leq m$, \textit{if} 
$g \in \text{H}^0(L^{\prime})$ \textit{vanishes at} $x_{i_1}, \ldots ,x_{i_l}$ 
\textit{then} $g = 0$. 

Choose a basis $f_1, \ldots ,f_m$ of $V$ such that $f_i(x_i) \neq 0$ and 
$f_i(x_j) = 0$ for $i \neq j$. Now, consider an element: 
\[
\xi = \underset{1\leq i_1<\cdots <i_r\leq m}{\textstyle \sum} 
g_{i_1\ldots i_r}\otimes f_{i_1}\wedge \ldots \wedge f_{i_r} \in 
\text{H}^0(L^{\prime})\otimes \overset{r}{\wedge} V 
\] 
such that $d(\xi ) = 0$. Recall that: 
\[
d(\xi ) = \underset{1\leq i_1<\cdots <i_{r-1}\leq m}{\textstyle \sum} 
({\textstyle \sum_{i=1}^m} f_ig_{ii_1\ldots i_{r-1}})\otimes 
f_{i_1}\wedge \ldots \wedge f_{i_{r-1}}, 
\]
using the alternate notation for $g_{ii_1\ldots i_{r-1}}$. Evaluating at $x_j$ 
the identity: 
\[
{\textstyle \sum_{i=1}^m}f_ig_{ii_1\ldots i_{r-1}} = 0,  
\]
one deduces that $g_{ji_1\ldots i_{r-1}}(x_j) = 0$. Consequently, for 
$1 \leq i_1 < \cdots < i_r \leq m$, $g_{i_1\ldots i_r}$ vanishes at 
$x_{i_1}, \ldots ,x_{i_r}$. It follows, from (ii) and from the hypothesis 
$r \geq l$, that $g_{i_1\ldots i_r} = 0$ hence $\xi = 0$.       
\end{proof} 
\vskip3mm 

\begin{proof}[Proof of Theorem 1] 
We shall prove the stability of $M_{d,V}$ by showing that, for $m$ large 
enough, Green's vanishing theorem (Lemma 1.2.) implies that $M_{d,V}$ 
satisfies the hypothesis of the Criterion of stability 1.1. More precisely, 
since $\text{rk}\, M_{d,V} = m-1$ it suffices to show that, for any two  
integers $0 < r < m-1$ and $a$: 
\[ 
\mu ((\overset{r}{\wedge} M_{d,V})\otimes {\mathcal O}_{\mathbb P}(a)) \leq 0 
\  \Rightarrow \  r \geq \text{h}^0({\mathcal O}_{\mathbb P}(a)) = P_n(a). 
\tag{1}
\] 
Now, $\mu ((\overset{r}{\wedge}M_{d,V})\otimes {\mathcal O}_{\mathbb P}(a)) = 
r\mu (M_{d,V}) + a = -\dfrac{rd}{m-1} + a$ hence: 
\[
\mu ((\overset{r}{\wedge} M_{d,V})\otimes {\mathcal O}_{\mathbb P}(a)) \leq 0 
\  \Leftrightarrow \  a \leq \frac{rd}{m-1}\  
\Leftrightarrow r \geq \frac{m-1}{d} a. 
\tag{2}  
\] 
From the inequality in the middle of (2) we deduce that $a < d$ hence 
$a \leq d-1$. We may also assume that $a \geq 1$. (1) would be true if one 
would show that: 
\[
\frac{m-1}{d} a > P_n(a) - 1,\  \text{for}\  1\leq a \leq d-1 
\tag{3} 
\]
which is equivalent to: 
\[
\frac{m-1}{d} > \frac{P_n(a) - 1}{a} = Q_{n-1}(a),\  \text{for} \   
1\leq a \leq d-1. 
\tag{4} 
\]

Now, $P_n(T)$ is a polynomial with positive coefficients hence $Q_{n-1}$ is a 
polynomial with positive coefficients, too. One deduces that the function 
$t \mapsto Q_{n-1}(t)$ is an increasing function for $t > 0$. Consequently, 
(4) is equivalent to: 
\[
\frac{m-1}{d} > Q_{n-1}(d-1). 
\tag{5} 
\] 
 Recalling the definition of $Q_{n-1}$, one deduces that (5)  
is equivalent to: 
\[
m > P_{n}(d-1) + Q_{n-1}(d-1). 
\tag{6} 
\]
Notice that, since $Q_{n-1}(d) > Q_{n-1}(d-1)$, (5) is satisfied by 
$m = P_n(d)$ hence, from (6), $P_n(d) > P_n(d-1) + Q_{n-1}(d-1)$. 

Similarly, if $m = P_n(d-1) + Q_{n-1}(d-1)$ then: 
\[
\mu ((\overset{r}{\wedge} M_{d,V})\otimes {\mathcal O}_{\mathbb P}(a)) < 0 
\  \Rightarrow \  r \geq \text{h}^0({\mathcal O}_{\mathbb P}(a)) 
\]
hence $M_{d,V}$ is semistable in this case.   
\end{proof} 
\vskip3mm 

{\bf 1.3. Example.}\quad 
Assume that $n = 2$, $d = 2$ and $m = 5$. Let $V$ be a b.p.f. 5-dimensional 
subspace of $S_2$ and let $M_V = M_{2,V}$ be the kernel of the evaluation 
morphism ${\mathcal O}_{{\mathbb P}^2}\otimes_k V\rightarrow 
{\mathcal O}_{{\mathbb P}^2}(2)$. Remark that, in this case, 
$m = P_2(1) + Q_1(1)$ hence, according to the last part of Theorem 1, 
$M_V$ is semistable. We shall prove the following two assertions: 

(a) \textit{For a general} $V$, $M_V$ \textit{is stable},\\ 
\hspace*{3mm} (b) $M_V$ \textit{does not satisfy the hypothesis of the 
Criterion of stability} 1.1. 

Indeed, $\mu (M_V) = -\dfrac{1}{2}$ hence $\mu (M_V\otimes 
{\mathcal O}_{\mathbb P}(a))$ is $\leq 0$ if and only if it is $< 0$, 
the same is true for $\mu (\overset{3}{\wedge} M_V\otimes 
{\mathcal O}_{\mathbb P}(a))$, and $\mu (\overset{2}{\wedge} M_V\otimes 
{\mathcal O}_{\mathbb P}(a)) \leq 0$ if and only if $a \leq 1$. It follows from 
the proof of Theorem 1 that if $r \in \{1,3\}$ and $\mu (\overset{r}{\wedge} 
M_V\otimes {\mathcal O}_{\mathbb P}(a)) \leq 0$ then 
$\text{H}^0(\overset{r}{\wedge} M_V\otimes {\mathcal O}_{\mathbb P}(a)) = 0$ 
and also that $\text{H}^0(\overset{2}{\wedge} M_V) = 0$. On the other hand, 
one deduces from the exact sequence: 
\[
0\longrightarrow (\overset{3}{\wedge} M_V)(-2) \longrightarrow 
{\mathcal O}_{\mathbb P}(-2)\otimes_k \overset{3}{\wedge} V \longrightarrow 
\overset{2}{\wedge} M_V \longrightarrow 0 
\]
that $\text{H}^0((\overset{2}{\wedge} M_V)(1)) \simeq 
\text{H}^1((\overset{3}{\wedge} M_V)(-1))$. Since $\overset{3}{\wedge} M_V 
\simeq M_V^{\ast}(-2)$ it follows, from Serre duality, that 
$\text{H}^1((\overset{3}{\wedge} M_V)(-1)) \simeq \text{H}^1(M_V)^{\ast}$. 
Using  the exact sequence: 
\[
0\longrightarrow M_V \longrightarrow {\mathcal O}_{{\mathbb P}^2}\otimes_k V 
\longrightarrow {\mathcal O}_{{\mathbb P}^2}(2) \longrightarrow 0 
\tag{1} 
\]
one deduces that $\text{H}^1(M_V) \simeq S_2/V$, hence it is 1-dimensional. 
This proves assertion (b). 

Now, if $M_V$ is not stable then it has a rank-2 coherent subsheaf 
$\mathcal F$, with $\text{det}\, {\mathcal F} \simeq 
{\mathcal O}_{\mathbb P}(-1)$ and such that $M_V/{\mathcal F}$ is torsion free. 
It follows, from the exact sequence $0 \rightarrow {\mathcal F} 
\rightarrow M_V \rightarrow M_V/{\mathcal F} \rightarrow 0$, that 
$\mathcal F$ is locally a 2-syzygy, hence locally free, and that the dual 
morphism $M_V^{\ast} \rightarrow {\mathcal F}^{\ast}$ is an epimorphism 
except at finitely many points. Since $M_V^{\ast}$ is a quotient of 
${\mathcal O}_{\mathbb P}\otimes_k V^{\ast}$, it follows that 
${\mathcal F}^{\ast} \simeq {\mathcal F}(1)$ is generated by global sections 
except at finitely many points. 
Since $\text{H}^0({\mathcal F}) = 0$ (because ${\mathcal F} \subset M_V$), 
$\mathcal F$ can be realized as an extension: 
\[
0 \longrightarrow {\mathcal O}_{\mathbb P}(-1) \longrightarrow {\mathcal F} 
\longrightarrow {\mathcal I}_{\Gamma} \longrightarrow 0 
\]
where $\Gamma$ is a 0-dimensional subscheme of ${\mathbb P}^2$. 
${\mathcal I}_{\Gamma}(1)$ must be generated by global sections except at 
finitely many points, hence $\text{h}^0({\mathcal I}_{\Gamma}(1)) \geq 2$, 
hence $\Gamma =$ 1 simple point, hence ${\mathcal F} \simeq 
{\Omega}_{\mathbb P}(1)$. 

Now, the inclusion ${\Omega}_{\mathbb P}(1) \rightarrow M_V$ determines (up to 
multiplication by a non-zero constant) an element $\xi \in \text{H}^1
(M_V(-1))$ (namely, a generator of the image of $\text{H}^1
({\Omega}_{\mathbb P}) \rightarrow \text{H}^1(M_V(-1))$) such that $h\xi = 0$, 
$\forall h \in S_1$. Applying $\text{Hom}_{{\mathcal O}_{\mathbb P}}(-,M_V)$ to 
the exact sequence: 
\[
0 \longrightarrow {\Omega}_{\mathbb P}(1) \longrightarrow 
{\mathcal O}_{\mathbb P}^3 \longrightarrow {\mathcal O}_{\mathbb P}(1) 
\longrightarrow 0 
\]
one deduces that $\xi = 0$ if and only if the morphism ${\Omega}_{\mathbb P}(1) 
\rightarrow M_V$ factorizes through ${\mathcal O}_{\mathbb P}^3$, which is 
not the case since $\text{H}^0(M_V) = 0$. Consequently, $\xi \neq 0$. 
One deduces, from the exact sequence (1), that $\xi$ corresponds to a non-zero 
element $f \in S_1$ such that $hf \in V$, $\forall h \in S_1$. We have thus 
proved that: 

(c) $M_V$ \textit{is not stable if and only if there exists} $0 \neq f \in 
S_1$ \textit{such that} $S_1f \subset V$. 

Now, let ${\mathbb P}(S_2)$ be the 5-dimensional projective space 
parametrizing the 1-dimensional quotients of $S_2$ (Grothendieck's convention) 
and let $v_2 : {\mathbb P}^2 \rightarrow {\mathbb P}(S_2)$ be the Veronese 
embedding. If $V \subset S_2$ is a 5-dimensional subspace then 
$S_2/V \in v_2({\mathbb P}^2)$ if and only if $V$ has a base point. If 
$0 \neq f \in S_1$  and if $L \subset {\mathbb P}^2$ is the line of equation 
$f = 0$ then the subset of ${\mathbb P}(S_2)$ consisting of the points 
$S_2/V$ with $V \supset S_1f$ is a 2-plane containing the conic $v_2(L)$, 
hence it is the linear span of that conic. Consequently, the subset of 
${\mathbb P}(S_2)$ consisting of the points $S_2/V$ with $V$ containing 
$S_1f$ for some non-zero $f \in S_1$ is exactly the {\it secant variety} 
$\text{Sec}\, v_2({\mathbb P}^2)$ of the Veronese embedding. As it is well 
known, $\text{Sec}\, v_2({\mathbb P}^2)$ is a determinantal cubic 
hypersurface in ${\mathbb P}(S_2)$. If $S_2/V \in {\mathbb P}(S_2) 
\setminus \text{Sec}\, v_2({\mathbb P}^2)$ then, by (c), $M_V$ is stable. 
This concludes the proof of assertion (a). We remark that if $S_2/V$ lies 
outside the secant variety of the Veronese embedding then the epimorphism 
${\mathcal O}_{{\mathbb P}^2}\otimes_k V \rightarrow {\mathcal O}_
{{\mathbb P}^2}(2)$ defines an embedding ${\varphi}_V : {\mathbb P}^2 
\rightarrow {\mathbb P}(V) \simeq {\mathbb P}^4$ and 
$M_V \simeq {\varphi}^{\ast}_V{\Omega}_{{\mathbb P}(V)}(1)$. 

We close the example by emphasizing another interesting property of $M_V$, 
namely: 

(d) \textit{If} $M_V$ \textit{is stable then it is uniform}. 

Indeed, it follows from (c) that, for every line $L \subset {\mathbb P}^2$,  
the kernel of the application $V \rightarrow \text{H}^0{\mathcal O}_L(2)$ 
has dimension $\leq 2$, hence the application is surjective. Restricting to 
$L$ the exact sequence (1) one deduces that $M_V\, \vert \, L \simeq 
{\mathcal O}_L^2 \oplus {\mathcal O}_L(-1)^2$. Computing the dimension of 
$\text{H}^j(M_V(-i))$ for $0 \leq i,j \leq 2$ one derives that the Beilinson 
monad of $M_V$ must have the form: 
\[
0 \longrightarrow {\Omega}_{\mathbb P}^2(2) \longrightarrow 
{}
{\Omega}_{\mathbb P}^1(1)^3 \longrightarrow {\mathcal O}_{\mathbb P} 
\longrightarrow 0
\]
hence $M_V(1)$ must be one of the bundles constructed by Elencwajg \cite{ele}. 
\vskip3mm     
     
\begin{proof}[Proof of Proposition 2] 
We may assume that $(n,d) \neq (2,2)$ since the remaining   
case is very easy. 
Let $Y \subset {\mathbb P}^n$ be any smooth complete intersection of $c$ 
hypersurfaces of degree $d$. We will show that, for any two integers 
$0 < r < \text{rk}\, N_{d,c}$ and $a$:  
\[
\mu ((\overset{r}{\wedge} N_{d,c})\otimes {\mathcal O}_Y(a)) \leq 0 \  
\Rightarrow \  \text{H}^0((\overset{r}{\wedge} N_{d,c})\otimes 
{\mathcal O}_Y(a)) = 0. 
\]
We have $\text{rk}\, N_{d,c} = \text{h}^0{\mathcal O}_Y(d) - 1 = 
P_n(d) - c - 1$ and $\text{det}\, N_{d,c} \simeq {\mathcal O}_Y(-d)$. It 
follows that: 
\[
\mu ((\overset{r}{\wedge} N_{d,c})\otimes {\mathcal O}_Y(a)) \leq 0 \  
\Leftrightarrow \  a \leq \frac{rd}{P_n(d)-c-1} \  
\Leftrightarrow \  r \geq \frac{P_n(d)-c-1}{d} a. 
\]
From the inequality in the middle one derives that $a \leq d-1$. In order to 
be able to apply Green's vanishing theorem (Lemma 1.2.) we would like to 
deduce from the last inequality that $r \geq \text{h}^0{\mathcal O}_Y(a) = 
P_n(a)$ (because $a \leq d-1$). We may assume that $a \geq 1$. As in the 
proof of Theorem 1, it would suffice to prove that: 
\[
\frac{P_n(d)-c-1}{d} > \frac{P_n(d-1)-1}{d-1} 
\]
which is equivalent to: 
\[
Q_{n-1}(d) - Q_{n-1}(d-1) > \frac{c}{d}\  . 
\]
But $Q_{n-1}(T)$ is a polynomial with positive coefficients, hence the same 
is true for the polynomial $Q_{n-1}(T+1) - Q_{n-1}(T)$. It follows that the  
function $t \mapsto Q_{n-1}(t+1) - Q_{n-1}(t)$ is an increasing function for 
$t > 0$ if $n \geq 3$ and a constant function if $n = 2$. One derives that: 
\[
Q_{n-1}(d) - Q_{n-1}(d-1) \geq Q_{n-1}(2) - Q_{n-1}(1) = 
\frac{n(n-1)}{2} \geq \frac{n-1}{2} \geq \frac{c}{d}\  . 
\]
Since $(n,d) \neq (2,2)$, at least one of the last two inequalities must be 
strict. 

Finally, as in the paper of Flenner \cite{fle}, using the fact that the 
relative Picard group of the universal family of complete intersections of 
$c$ hypersurfaces of degree $d$ in ${\mathbb P}^n$ is generated by the 
pullback of ${\mathcal O}_{{\mathbb P}^n}(1)$ and considering the relative 
Harder-Narasimhan filtration of the vector bundle on this universal family 
patching together the bundles $N_{d,c}$, one deduces that, for a {\it general}  
$Y$, $N_{d,c}$ is semistable with respect to ${\mathcal O}_Y(1)$. 
\end{proof} 
\vskip3mm 

{\bf 1.4. Question.}\quad 
Assume that $\text{char}\, k = p > 0$. Is it true that the bundle 
$N_{d,c}$ from Proposition 2 is {\it strongly semistable}, in the sense that 
all of its iterated Frobenius pullbacks are semistable? A positive answer to 
this question would lead to an improvement of the estimate in a recent 
restriction theorem of Langer \cite{lan}, Theorem 2.1., in the same way the 
Proposition 2 leads to an improvement of the estimate in Flenner's 
restriction theorem. 

\section{Brenner's criterion of stability}

Let $u_1,\ldots ,u_m \in S = k[X_0,\ldots ,X_n]$ be distinct monomials of 
degrees $d_1 \leq \cdots \leq d_m$ with no common factor. 
$V := ku_1 + \dots + ku_m$ is a {\it graded} $k$-vector subspace of $S$. 
Let $(K_{\bullet},{\delta}_K)$ be the Koszul complex of graded $S$-modules 
defined by $u_1, \ldots ,u_m$: 
\[
0 \rightarrow S\otimes_k \overset{m}{\wedge} V \rightarrow \cdots 
\rightarrow S\otimes_k V \rightarrow S \rightarrow 0 .
\]
Notice that, as a graded $S$-module, $S\otimes_k V \simeq 
\bigoplus_{i=1}^mS(-d_i)$. Let $\varphi : \bigoplus_{i=1}^m
{\mathcal O}_{{\mathbb P}^n}(-d_i) \rightarrow {\mathcal O}_{{\mathbb P}^n}$ be 
the morphism of ${\mathcal O}_{{\mathbb P}^n}$-modules defined by 
$u_1, \ldots ,u_m$, i.e., the sheafification of the morphism of graded 
$S$-modules $S\otimes_k V \rightarrow S$, and let $\mathcal F$ be the kernel 
of $\varphi $. For $I \subseteq \{1,\ldots ,m\}$ with $\text{card}\, I 
\geq 2$, let ${\mathcal F}_I$ denote the kernel of $\varphi \, \vert \, 
\bigoplus_{i\in I}{\mathcal O}_{{\mathbb P}^n}(-d_i) \rightarrow 
{\mathcal O}_{{\mathbb P}^n}$. Brenner's criterion of stability is the 
following: 

\vskip3mm 
{\bf 2.1. Theorem.} (Brenner \cite{bre}, Theorem 6.3.)\quad 
\textit{Using the above hypotheses and notation, if} ${\mathcal F}^{\prime}$ 
\textit{is a coherent subsheaf of} $\mathcal F$ \textit{of rank} $r$ 
\textit{then}: 
\[
\mu ({\mathcal F}^{\prime}) \leq \text{max}\, \{\mu ({\mathcal F}_I)\, \vert 
\, I\subseteq \{1,\ldots ,m\},\  2\leq \text{card}\, I \leq r+1 \}. 
\]
\vskip3mm 

Before giving a proof of this theorem, we introduce some more notation. 
Let $W := k^m$, let $e_1, \ldots ,e_m$ be the canonical $k$-basis of $W$ and 
let $(C_{\bullet},{\delta}_C)$ be the Koszul complex $0 \rightarrow 
\overset{m}{\wedge} W \rightarrow \cdots \rightarrow W \rightarrow k 
\rightarrow 0$ defined by ${\delta}_C(e_i) = 1$, $i = 1, \ldots ,m$. For 
$I \subseteq \{1,\ldots ,m\}$ let $e_I$ denote the {\it exterior monomial} 
$e_{i_1} \wedge \ldots \wedge e_{i_r}$, where $I = \{i_1, \ldots ,i_r\}$ and 
$i_1 < \dots < i_r$. For $\omega = \sum_{i_1<\dots <i_r}c_{i_1\ldots i_r}
e_{i_1}\wedge \ldots \wedge e_{i_r} \in \overset{r}{\wedge} W$ we define: 
\[
\text{Supp}(\omega) := \{e_{i_1}\wedge \ldots \wedge e_{i_r}\, \vert \, 
i_1<\dots <i_r\  \text{and}\  c_{i_1\ldots i_r}\neq 0\}.
\]
We also consider the Koszul complex of graded $S$-modules $(L_{\bullet},
{\delta}_L) := S\otimes_k (C_{\bullet},{\delta}_C)$. The morphism of graded 
$S$-modules $f: S\otimes_k V \rightarrow S\otimes_k W$, 
$f(1\otimes u_i) := u_i\otimes e_i$, $i = 1,\ldots ,m$, extends to a morphism 
of complexes $\bigwedge f : K_{\bullet} \rightarrow L_{\bullet}$. 

Now, a key point in Brenner's proof of his criterion is the following 
result (compare with \cite{bre}, Lemma 6.2., and the last part of the proof 
of \cite{bre}, Theorem 6.3.): 

\vskip3mm 
{\bf 2.2. Lemma.} \textit{Let} $\omega$ \textit{be a non-zero decomposable 
element of} $\overset{r}{\wedge} W$ \textit{such that} ${\delta}_C(\omega) 
= 0$. \textit{Then there exist mutually disjoint subsets} $I_1, \ldots , 
I_s$ \textit{of} $\{1, \ldots ,m\}$, \textit{with} $\text{card}\, I_i 
\geq 2$, $i = 1, \ldots ,s$, \textit{and} $\text{card}\, I_1 + \dots 
+ \text{card}\, I_s = r+s$, \textit{such that}: 
\[
\text{Supp}(\omega) \supseteq \text{Supp}({\delta}_C(e_{I_1})\wedge \ldots 
\wedge {\delta}_C(e_{I_s})) .
\]
\vskip3mm 

\begin{proof}
Assume that $\omega = w_1\wedge \ldots \wedge w_r$ and let $W^{\prime}$ be the 
subspace of $W$ spanned by $w_1, \ldots ,w_r$. Using Gaussian elimination, 
one can find a subset $I = \{i_1, \ldots ,i_r\}$, $i_1 < \dots < i_r$ of 
$\{1, \ldots ,m\}$ and a $k$-basis $w_1^{\prime}, \ldots ,w_r^{\prime}$ of 
$W^{\prime}$ such that, denoting by $I^{\prime}$ the complement 
$\{1, \ldots ,m\} \setminus I$ of $I$ and by $I_p^{\prime}$ the set 
$\{j\in I^{\prime}\, \vert \, j>i_p\}$, $p = 1, \ldots ,r$, one has: 
\[
w_p^{\prime} = e_{i_p} + {\textstyle \sum_{j\in I_p^{\prime}}}a_{pj}e_j,\  
p=1,\ldots ,r.
\]
${\omega}^{\prime} := w_1^{\prime}\wedge \ldots \wedge w_r^{\prime}$ differs 
from $\omega$ by a non-zero multiplicative constant. Since: 
\[
0 = {\delta}_C({\omega}^{\prime}) = {\textstyle \sum}(-1)^{p-1}
{\delta}_C(w_p^{\prime})w_1^{\prime}\wedge \ldots \wedge  
{\widehat {w_p^{\prime}}} \wedge \ldots \wedge w_r^{\prime} 
\]
it follows that, for $1\leq p\leq r$,  ${\delta}_C(w_p^{\prime}) = 0$ hence 
$\exists j \in I_p^{\prime}$ such that $a_{pj} \neq 0$. Let $j_p$ be the least 
$j$ with this property and let $a_p := a_{pj_p}$, $p = 1, \ldots ,r$. 

{\bf Claim:} $\eta := (e_{i_1}+a_1e_{j_1})\wedge \ldots \wedge 
(e_{i_r}+a_re_{j_r})$ \textit{and} ${\omega}^{\prime} - \eta$ \textit{have 
disjoint supports}. 

Indeed, consider an exterior monomial $\varepsilon = e_{l_1}\wedge \ldots 
\wedge e_{l_r}$, $l_1 < \dots < l_r$. We associate to $\varepsilon$ the subset 
$A := \{1\leq p \leq r\, \vert \, i_p\in \{l_1,\ldots ,l_r\}\}$ of 
$\{1, \ldots ,r\}$ and its complement $A^{\prime} := \{1, \ldots ,r\} 
\setminus A$. 

If $\varepsilon \in \text{Supp}(\eta)$ then there exists a permutation 
$\sigma \in \mathfrak{S}_r$ such that $l_{\sigma (p)} \in \{i_p,j_p\}$, 
$p = 1, \ldots ,r$. It follows that $A = \{p\, \vert \, l_{\sigma (p)} = 
i_p\}$ and $A^{\prime} = \{p\, \vert \, l_{\sigma (p)} = j_p\}$ hence 
$l_1 + \dots + l_r = l_{\sigma (1)} + \dots + l_{\sigma (r)} = 
\sum_{p\in A} i_p + \sum_{p\in A^{\prime}} j_p$. 

On the other hand, if $\varepsilon \in \text{Supp}({\omega}^{\prime} - \eta)$ 
then there exists a permutation $\tau \in \mathfrak{S}_r$ such that 
$l_{\tau (p)} \in \{i_p\} \cup \{j\in I^{\prime}\, \vert \, j \geq j_p\}$, 
$p = 1, \ldots ,r$, and, for at least one $p$, $l_{\tau (p)} \in \{j\in 
I^{\prime}\, \vert \, j > j_p\}$. It follows that $A = \{p\, \vert \, 
l_{\tau (p)} = i_p\}$ and $A^{\prime} = \{p\, \vert \, l_{\tau (p)} \in 
I^{\prime},\  l_{\tau (p)} \geq j_p\}$ hence $l_1 + \dots + l_r = 
l_{\tau (1)} + \ldots + l_{\tau (r)} > \sum_{p\in A}i_p + \sum_{p\in A^{\prime}}
j_p$. The Claim is proved. 

It follows from the Claim that $\text{Supp}({\omega}^{\prime}) \supseteq 
\text{Supp}(\eta) = \text{Supp}((e_{i_1}-e_{j_1}) \wedge \ldots \wedge 
(e_{i_r}-e_{j_r}))$. Let $j_1^{\prime}, \ldots ,j_s^{\prime}$ be the distinct 
elements of the set $\{j_1, \ldots ,j_r\}$ and, for $1 \leq q \leq s$, let 
$I_q := \{j_q^{\prime}\}\cup \{i_p\, \vert \, 1 \leq p \leq r,j_p = 
j_q^{\prime}\}$. Using the relation: 
\[
(e_{l_1}-e_{l_0})\wedge \ldots \wedge (e_{l_p}-e_{l_0}) = 
{\delta}_C(e_{l_0}\wedge e_{l_1}\wedge \ldots \wedge e_{l_p})
\]
one deduces that: 
$
(e_{i_1} -e_{j_1})\wedge \ldots \wedge (e_{i_r}-e_{j_r}) = 
\pm {\delta}_C(e_{I_1})\wedge \ldots \wedge {\delta}_C(e_{I_s}).
$
\end{proof}
\vskip3mm 

{\bf 2.3. Remark.} \textit{The conclusion of Lemma} 2.2. 
\textit{is no more valid if one 
assumes only that} ${\delta}_C(\omega) = 0$: \textit{the element}
$
\omega = e_1\wedge e_2 + e_2\wedge e_3 + e_3\wedge e_4 + e_4\wedge e_5 + 
e_5\wedge e_1 \in \overset{2}{\wedge} V
$
\textit{provides a counterexample}. 
\vskip3mm 

In order to check whether an element $\omega \in \overset{r}{\wedge} W$ is 
{\it decomposable} or not one can use the following well known criterion: 
          
\vskip3mm 
{\bf 2.4. Lemma.}\quad 
\textit{An element} $\omega = \sum_{i_1<\dots <i_r}c_{i_1\ldots i_r}
e_{i_1}\wedge \ldots \wedge e_{i_r} \in \overset{r}{\wedge} W$ 
\textit{is decomposable} (\textit{i.e.}, \textit{of the form} 
$w_1\wedge \ldots \wedge w_r$ \textit{for some vectors} $w_1, \ldots ,w_r 
\in W$) \textit{if and only if}: 
\[
{\textstyle \sum_{p=0}^r}(-1)^pc_{i_1\ldots i_{r-1}j_p}\cdot 
c_{j_0\ldots {\widehat {j_p}}\ldots j_r} = 0,\  \forall \, i_1<\dots <i_{r-1},\  
\forall \, j_0<\dots <j_r. 
\]
\vskip3mm 

\begin{proof} 
One considers the {\it contraction pairing}: 
\[
\langle -,- \rangle : \overset{r}{\wedge} W\otimes \overset{r-1}{\wedge} 
W^{\ast} \longrightarrow W\otimes \overset{r-1}{\wedge} W\otimes 
\overset{r-1}{\wedge} W^{\ast} \longrightarrow W. 
\]
Let $W^{\prime}$ be the $k$-vector subspace of $W$ spanned by the vectors 
$\langle \omega ,\psi \rangle$, $\psi \in \overset{r-1}{\wedge} W^{\ast}$. 

(i) \textit{We assert, firstly, that} 
$\omega \in \overset{r}{\wedge} W^{\prime}$. 

Indeed, let $w_1,\ldots ,w_m$ be a $k$-basis of $W$ such that $w_1,\ldots ,
w_q$ is a $k$-basis of $W^{\prime}$. Let $w_1^{\ast}, \ldots ,w_m^{\ast}$ be 
the dual basis of $W^{\ast}$. If we express $\omega$ as 
$\sum_{i_1<\dots <i_r}b_{i_1\ldots i_r}w_{i_1}\wedge \ldots \wedge w_{i_r}$ 
then: 
\[
\langle \omega,w_{j_1}^{\ast}\wedge \ldots \wedge w_{j_{r-1}}^{\ast}\rangle = 
{\textstyle \sum_{j=1}^m}b_{jj_1\ldots j_{r-1}}w_j = 
(-1)^{r-1}{\textstyle \sum_{j=1}^m}b_{j_1\ldots j_{r-1}j}w_j. 
\]
If $\omega \notin \overset{r}{\wedge} W^{\prime}$ then there exist 
$i_1< \dots < i_{r}$ with $i_{r} > q$ such that $b_{i_1\ldots i_{r}} 
\neq 0$. In this case, $\langle \omega,w_{i_1}^{\ast}\wedge \ldots \wedge 
w_{i_{r-1}}^{\ast}\rangle \notin W^{\prime}$, a contradiction. 

(ii) \textit{Next}, \textit{we assert that} $\omega$ \textit{is decomposable 
if and only if} $\omega \wedge w^{\prime} = 0$, $\forall \, w^{\prime} \in 
W^{\prime}$. 

Indeed, if $\omega = w_1\wedge \ldots \wedge w_r$ then, extending $w_1, \ldots 
,w_r$ to a $k$-basis of $W$ and considering the dual basis of $W^{\ast}$, 
one deduces that $W^{\prime}$ is generated by $w_1, \ldots ,w_r$, hence 
$\omega \wedge w^{\prime} = 0$, $\forall \, w^{\prime} \in W^{\prime}$. 

Conversely, assume that $\omega \wedge w^{\prime} = 0$, $\forall \, w^{\prime} 
\in W^{\prime}$. In order to show that $\omega$ is decomposable it suffices, 
by (i), to show that $\text{dim}_kW^{\prime} \leq r$. If $q := \text{dim}_k 
W^{\prime} > r$ then, as the pairing $\overset{r}{\wedge} W^{\prime} 
\otimes \overset{q-r}{\wedge} W^{\prime} \rightarrow \overset{q}{\wedge} 
W^{\prime}$ is non-degenerate, one gets that $\omega = 0$, a contradiction. 

(iii) \textit{Finally}, if $e_1^{\ast}, \ldots ,e_m^{\ast}\in W^{\ast}$  
is the dual of the canonical basis $e_1, \dots ,e_m$ of $W$ then the system 
of equations: 
\[
\omega \wedge \langle \omega ,e_{i_1}^{\ast}\wedge \ldots \wedge 
e_{i_{r-1}}^{\ast}\rangle = 0,\  i_1 < \dots < i_{r-1}, 
\]
is equivalent to the system of equations from the statement.    
\end{proof}
\vskip3mm 

\begin{proof}[Proof of Theorem 2.1] 
We may assume that ${\mathcal F}/{\mathcal F}^{\prime}$ is torsion free. One 
has $\text{deg}\, {\mathcal F}^{\prime} = -d$ for some $d \in {\mathbb N}$. 
Consider an open subset $U$ of ${\mathbb P}^n$, with $\text{codim}
({\mathbb P}^n\setminus U,{\mathbb P}^n) \geq 2$, such that $\varphi \, \vert 
\, U$ is an epimorphism and $({\mathcal F}/{\mathcal F}^{\prime})\, \vert \, 
U$ is locally free. In this case, ${\mathcal F}\, \vert \, U$ is locally 
free, one has an exact sequence: 
\[
0 \longrightarrow \overset{r}{\textstyle \bigwedge} 
({\mathcal F}\, \vert \, U) \longrightarrow 
\overset{r}{\textstyle \bigwedge} 
({\textstyle \bigoplus_{i=1}^m}{\mathcal O}_U(-d_i)) 
\longrightarrow \overset{r-1}{\textstyle \bigwedge} 
({\textstyle \bigoplus_{i=1}^m}{\mathcal O}_U(-d_i)) 
\]
and ${\mathcal F}^{\prime}\, \vert \, U$ is locally a direct summand of 
${\mathcal F}\, \vert \, U$, hence locally free. Moreover, 
$\overset{r}{\bigwedge} ({\mathcal F}^{\prime}\, \vert \, U) \simeq 
{\mathcal O}_U(-d)$. 

Considering the Koszul complex $(K_{\bullet},{\delta}_K)$ associated to 
the monomials $u_1, \ldots ,u_m$ (see the beginning of the section), the 
inclusion $\overset{r}{\bigwedge} ({\mathcal F}^{\prime}\, \vert \, U) 
\hookrightarrow \overset{r}{\bigwedge} ({\mathcal F}\, \vert \, U)$ defines 
a homogeneous element $\xi$ of $K_r$, of degree $d$, such that 
${\delta}_K(\xi) = 0$. One can write: 
\[
\xi = \underset{i_1<\dots <i_r}{\textstyle \sum} P_{i_1\ldots i_r}\otimes 
u_{i_1}\wedge \ldots \wedge u_{i_r} 
\]
with $P_{i_1\ldots i_r} \in S$ homogeneous polynomial of degree 
$d - d_{i_1} - \cdots - d_{i_r}$.  The morphisms $u_i \cdot - : 
{\mathcal O}_{{\mathbb P}^n}(-d_i) \rightarrow  
{\mathcal O}_{{\mathbb P}^n}$, $i = 1, \ldots ,m$, define a morphism 
$\bigoplus_{i=1}^m{\mathcal O}_{{\mathbb P}^n}(-d_i) \rightarrow  
{\mathcal O}_{{\mathbb P}^n}^m$ which is just  
the sheafification of the morphism of graded $S$-modules 
$f : S\otimes_k V \rightarrow S\otimes_k W$, $f(1\otimes u_i) = u_i\otimes 
e_i$, $i = 1, \ldots ,m$, considered after the statement of Theorem 2.1. 
It is an isomorphism on the open set $U_{0\ldots n} := {\mathbb P}^n 
\setminus \bigcup_{i=0}^n\{X_i = 0\}$. Let: 
\[
\eta := (\overset{r}{\textstyle \bigwedge} f)(\xi) = 
\underset{i_1<\dots <i_r}{\textstyle \sum} u_{i_1}\ldots u_{i_r}P_{i_1\ldots i_r} 
\otimes e_{i_1}\wedge \ldots \wedge e_{i_r} \in L_r\  .
\]
One has ${\delta}_L(\eta) = (\overset{r-1}{\bigwedge} f)({\delta}_K(\xi)) 
= 0$. Let us denote $u_{i_1}\ldots u_{i_r}P_{i_1\ldots i_r}$ by 
$Q_{i_1\ldots i_r}$. It is a homogeneous polynomial of degree $d$. It follows 
from  the definition of $\xi$ that, $\forall \, x \in U$, 
$\eta (x) = \sum_{i_1<\dots <i_r}Q_{i_1\ldots i_r}(x)e_{i_1}\wedge \ldots \wedge 
e_{i_r}$ is a {\it decomposable} element of $\overset{r}{\wedge} W$. One 
deduces from Lemma 2.4. that: 
\[
{\textstyle \sum_{p=0}^r}(-1)^pQ_{i_1\ldots i_{r-1}j_p}(x)\cdot  
Q_{j_0\ldots \widehat{j_p}\ldots j_r}(x) = 0,\  
\forall \, i_1<\dots <i_{r-1},\  \forall \, j_0<\dots <j_r, 
\]
hence, for the same multiindices: 
\[
{\textstyle \sum_{p=0}^r}(-1)^pQ_{i_1\ldots i_{r-1}j_p}\cdot 
Q_{j_0\ldots \widehat{j_p}\ldots j_r} = 0\  \text{in}\  S_{2d}\  . 
\tag{*}
\]

Consider, now, a {\it monomial order} on $S$, let's say the lexicographic 
order with $X_0 > \dots > X_n$. Let $u$ be the largest monomial among the 
initial monomials of the polynomials $Q_{i_1\ldots i_r}$ which are non-zero. 
One can write: 
\[
Q_{i_1\ldots i_r} = c_{i_1\ldots i_r}u + \text{a linear combination of monomials} 
<_{\text{lex}}\, u 
\]
where, of course, $c_{i_1\ldots i_r} = 0$ if $Q_{i_1\ldots i_r} = 0$ or if its 
initial monomial is $<_{\text{lex}}\, u$. Consider the element: 
\[
\omega := \underset{i_1<\dots <i_r}{\textstyle \sum}c_{i_1\ldots i_r} 
e_{i_1}\wedge \ldots \wedge e_{i_r} \in \overset{r}{\wedge} W. 
\]
Interpreting $\eta$ as a homogeneous polynomial of degree $d$ with 
coefficients in $\overset{r}{\wedge} W$, $\omega$ is just the (exterior) 
coefficient of $u$ in $\eta$. Since ${\delta}_L(\eta) = 0$, it follows that 
${\delta}_C(\omega) = 0$. Moreover, evaluating the coefficient of $u^2$ in 
the relations (*), one deduces that $\omega$ satisfies the relations from 
Lemma 2.4. hence it is a {\it decomposable} element of $\overset{r}{\wedge} 
W$. 

Finally, according to Lemma 2.2., there exist mutually disjoint subsets 
$I_1, \ldots ,I_s$ of $\{1,  \ldots ,m\}$, with $\text{card}\, I_i \geq 2$, 
$i = 1, \ldots ,s$, and with $\text{card}\, I_1 + \dots + \text{card}\, I_s 
= r+s$, such that: 
\[
\text{Supp}(\omega) \supseteq \text{Supp}({\delta}_C(e_{I_1})\wedge \ldots 
\wedge {\delta}_C(e_{I_s})). 
\]
We remark that if $e_{i_1}\wedge \ldots \wedge e_{i_r} \in \text{Supp}
(\omega)$ then $u_{i_1}\ldots u_{i_r}$ divides $u$. Let $I = I_1 \cup \ldots 
\cup I_s$. The products of monomials corresponding to the elements of 
$\text{Supp}({\delta}_C(e_{I_1})\wedge \ldots \wedge {\delta}_C(e_{I_s}))$ are: 
\[
\frac{{\textstyle \prod_{i\in I}}u_i}{u_{i_1}\ldots u_{i_s}}\  \text{with}\  
i_1\in I_1, \ldots ,i_s\in I_s\  . 
\]
Their least common multiple is: 
\[
\frac{{\textstyle \prod_{i\in I}}u_i}{\text{gcd}\, \{u_{i_1}\ldots u_{i_s}\, 
\vert \, i_1\in I_1, \ldots , i_s\in I_s\}}\  . 
\]
But $\text{gcd}\, \{u_{i_1}\ldots u_{i_s}\, \vert \, i_1\in I_1, \ldots ,i_s 
\in I_s\} = \text{gcd}\, \{u_{i_1}\, \vert \, i_1\in I_1\}\cdot \ldots \cdot 
\text{gcd}\, \{u_{i_s}\, \vert \, i_s\in I_s\}$. One deduces that $u$ is 
divisible by the product: 
\[
\frac{{\textstyle \prod_{i_1\in I_1}}u_{i_1}}{\text{gcd}\, \{u_{i_1}\, \vert 
\, i_1\in I_1\}}\cdot \ldots \cdot 
\frac{{\textstyle \prod_{i_s\in I_s}}u_{i_s}}{\text{gcd}\, \{u_{i_s}\, \vert 
\, i_s\in I_s\}}\  .
\]
The degree of the $p$th factor of this product equals $-\text{deg}\, 
{\mathcal F}_{I_p}$. Let us denote by $\mu$ the right hand side of the 
inequality from the statement of the theorem. Then $-\text{deg}\, 
{\mathcal F}_{I_p} \geq -(\text{card}\, I_p - 1)\mu$. It follows that: 
\[
d = \text{deg}(u) \geq -((\text{card}\, I_1 -1) + \dots + (\text{card}\, I_s 
- 1))\mu = -r\mu 
\]
and, consequently, $\mu ({\mathcal F}^{\prime}) := -\dfrac{d}{r} 
\leq \mu$.     
\end{proof} 
\vskip3mm 

One gets immediately the following: 

\vskip3mm 
{\bf 2.5. Corollary.} (Brenner \cite{bre}, Corollary 6.4.)\quad 
\textit{The sheaf} $\mathcal F$ \textit{defined at the beginning of this 
section is} (\textit{semi})\textit{stable if and only if}, \textit{for 
any integers} $2 \leq r \leq m-1$ \textit{and} $1 \leq i_1 < \dots < i_r 
\leq m$, \textit{one has}: 
\[
\frac{d_{i_1\ldots i_r}-d_{i_1}-\dots -d_{i_r}}{r-1}\  (\leq) < 
\frac{-d_1-\dots -d_m}{m-1}
\]
\textit{where} $d_{i_1\ldots i_r}$ \textit{is the degree of the greatest 
common divisor of} $u_{i_1},\ldots ,u_{i_r}$. 
\vskip3mm 

In the particular case $d_1=\dots =d_m=d$ one gets: 

\vskip3mm 
{\bf 2.6. Corollary.} (Brenner \cite{bre}, Corollary 6.5.)\quad 
\textit{If} $d_1=\dots =d_m=d$ \textit{then}, \textit{recalling that} 
$V := ku_1 + \dots + ku_m \subseteq S_d$, \textit{the sheaf} $\mathcal F$ 
\textit{is} (\textit{semi})\textit{stable if and only if} $\forall \, 
1 \leq e \leq d-1$, $\forall \, u \in S_e$ \textit{monomial}: 
\[
\frac{\text{dim}(V:u)-1}{d-e}\  (\leq) < \frac{\text{dim}\, V-1}{d}, 
\]
\textit{where} $(V:u) := \{f\in S_{d-e}\, \vert \, uf\in V\}$. 
\vskip3mm 

\begin{proof} 
The inequalities from Corollary 2.5. become, in this particular case: 
\[
\frac{d_{i_1\ldots i_r}-rd}{r-1}\  (\leq) < \frac{-md}{m-1} 
\]
which is equivalent to: 
\[
\frac{r-1}{d-d_{i_1\ldots i_r}}\  (\leq) < \frac{m-1}{d}\  . 
\]
Now, if $u = \text{gcd}(u_{i_1},\ldots ,u_{i_r})$ then $ku_{i_1} + \dots + 
ku_{i_r} \subseteq u(V:u)$, hence $r \leq \text{dim}(V:u)$. On the other 
hand, if $u$ is a monomial of degree $e$ and if $u(V:u) = ku_{i_1} + \dots + 
ku_{i_r}$ then $d_{i_1\ldots i_r} \geq e$. 
\end{proof} 
\vskip3mm

\section{Applications of Brenner's criterion} 

Throughout this section, we shall denote by $S^{\prime}$ the  
subalgebra $k[X_0,\ldots ,X_{n-1}]$ of the polynomial algebra 
$S = k[X_0,\ldots ,X_n]$. 

\begin{proof}[Proof of Theorem 4] 
We will show, by induction on $n\geq 2$, that for any integers $d \geq 1$ 
and $n+1 \leq m \leq P_n(d)$ there exists an $m$-dimensional b.p.f. monomial 
subspace $V$ of $S_d$ satifying the strict inequalities from Corollary 2.6., 
except for $n=2$, $d=2$, $m=5$, when $V$ satisfies only the non-strict 
inequalities. The case $n = 2$ is the main result of the paper of 
Costa, Macias Marques and Mir\'{o}-Roig \cite{cmm} recalled in Theorem 3 
from the Introduction (see, also, the remarks following this proof). 

For the proof of the {\it induction step} $(n-1)\rightarrow n\geq 3$ we use 
induction on $d$. In the case $d=1$ we have nothing to check. 
We shall divide the proof of the {\it induction step} $(d-1)\rightarrow 
d\geq 2$ into three cases. 

{\bf Case 1:} $n+1 \leq m \leq P_{n-1}(d)+1$. 

By the induction hypothesis on $n$, there exists an $(m-1)$-dimensional 
b.p.f. monomial subspace $V^{\prime}$ of $S^{\prime}_d$ such that, 
$\forall \, 1 \leq e \leq d-1$, $\forall \, u^{\prime} \in S^{\prime}_e$ 
monomial, one has:
\[
\frac{\text{dim}(V^{\prime}:u^{\prime})-1}{d-e} \leq 
\frac{\text{dim}\, V^{\prime}-1}{d}\  .
\]
Notice that we require only that $V^{\prime}$ satisfy the non-strict 
inequalities from Corollary 2.6. We take $V := V^{\prime} + kX_n^d$. Consider, 
now, an integer $1 \leq e \leq d-1$ and a monomial $u \in S_e$. 

(i) If $u \in S^{\prime}_e$ then $(V:u) = (V^{\prime}:u)$ hence: 
\[
\frac{\text{dim}(V:u)-1}{d-e} \leq \frac{\text{dim}\, V^{\prime}-1}{d} 
< \frac{\text{dim}\, V-1}{d}\  . 
\]
\hspace*{3mm} (ii) If $u=X_n^e$ then $(V:u)=kX_n^{d-e}$, hence 
$\text{dim}(V:u)-1 = 0$. 

(iii) If $u\in S_{e-1}X_n\setminus \{X_n^e\}$ then $(V:u) = (0)$. 

{\bf Case 2:} $P_{n-1}(d)+2 \leq m \leq P_{n-1}(d)+n+1$. 

In this case, $m=P_{n-1}(d)+n+1-l$ for some $0 \leq l \leq n-1$. We take 
$V := S^{\prime}_d + (kX_l^{d-1} + kX_{l+1}^{d-1} + \dots + kX_n^{d-1})X_n$. 
Consider an integer $1 \leq e \leq d-1$ and a monomial $u \in S_e$. 

(i) If $u \in S^{\prime}_e$ then $\text{dim}(V:u) \leq \text{dim}
(S^{\prime}_d:u) + 1 = \text{dim}\, S^{\prime}_{d-e} + 1$ hence: 
\[
\frac{\text{dim}(V:u)-1}{d-e} \leq \frac{P_{n-1}(d-e)}{d-e}\  . 
\]
On the other hand: 
\[
\frac{\text{dim}\, V-1}{d}\geq \frac{P_{n-1}(d)+1}{d} > 
\frac{P_{n-1}(d)}{d} 
\]
hence it suffices to show that: 
\[
\frac{P_{n-1}(d-e)}{d-e} \leq \frac{P_{n-1}(d)}{d}\  .
\]
In order to prove this inequality it suffices to show that: 
\[
\frac{P_{n-1}(\delta -1)}{\delta -1} \leq \frac{P_{n-1}(\delta)}{\delta}, 
\  \text{for}\  2\leq \delta \leq d. 
\]
One checks easily that the last inequality is equivalent to: 
\[
\delta \geq \frac{n-1}{n-2} = 1 + \frac{1}{n-2} 
\]
which is true because $n \geq 3$ and $\delta \geq 2$. 

(ii) If $u \in S_{e-1}X_n$ and $e \geq 2$ then $\text{dim}(V:u) - 1 \leq 0$. 

(iii) If $u=X_n$ then $(V:u) = kX_l^{d-1} + \dots + kX_n^{d-1}$ hence: 
\[
\frac{\text{dim}(V:u)-1}{d-1} = \frac{n-l}{d-1} < 
\frac{P_{n-1}(d) + n - l}{d} = \frac{\text{dim}\, V-1}{d}
\]
because $n-l \leq n < \dfrac{n(n+1)}{2} = P_{n-1}(2) \leq (d-1)P_{n-1}(d)$. 

{\bf Case 3:} $P_{n-1}(d)+n+1 \leq m \leq P_n(d)$. 

Since $P_n(d) - P_{n-1}(d) = P_n(d-1)$, we have $m = P_{n-1}(d) + l$, with 
$n+1 \leq l \leq P_n(d-1)$. By the induction hypothesis on $d$, there exists 
an $l$-dimensional b.p.f. monomial subspace $W$ of $S_{d-1}$ such that 
$\forall \, 1 \leq e \leq d-2$, $\forall \, u \in S_e$ monomial: 
\[
\frac{\text{dim}(W:u)-1}{d-1-e} \leq \frac{\text{dim}\, W-1}{d-1}\  . 
\]
Notice, again, that we require only that $W$ satisfy the non-strict 
inequalities from Corollary 2.6. We take $V := S^{\prime}_d + WX_n$. Consider 
an integer $1 \leq e \leq d-1$ and a monomial $u \in S_e$. 

(i) If $u \in S^{\prime}_e$ then $(V:u) = S^{\prime}_{d-e} + (W:u)X_n$, hence 
we have to show that: 
\[
\frac{P_{n-1}(d-e)+\text{dim}(W:u)-1}{d-e} < 
\frac{P_{n-1}(d)+\text{dim}\, W-1}{d}\  .
\tag{*} 
\]
We have already proved that: 
\[
\frac{P_{n-1}(d-e)}{d-e} \leq \frac{P_{n-1}(d)}{d}\  .
\]
Now, if $e = d-1$ then $\text{dim}(W:u) - 1 \leq 0$ and $\text{dim}\, 
W - 1 \geq n$ hence (*) is fulfilled. If $1 \leq e \leq d-2$ then: 
\[
\text{dim}(W:u)-1 \leq \frac{d-1-e}{d-1}(\text{dim}\, W-1) < 
\frac{d-e}{d}(\text{dim}\, W-1) 
\]
hence: 
\[
\frac{\text{dim}(W:u)-1}{d-e} < \frac{\text{dim}\, W-1}{d} 
\]
and (*) is again fulfilled. 

(ii) If $u\in S_{e-1}X_n$ and $e \geq 2$ then $(V:u) = (W:(u/X_n))$ hence: 
\[
\frac{\text{dim}(V:u)-1}{d-e} = 
\frac{\text{dim}(W:(u/X_n))-1}{d-1-(e-1)} \leq 
\frac{\text{dim}\, W-1}{d-1}\  . 
\]
We would like to show that: 
\[
\frac{\text{dim}\, W-1}{d-1} < \frac{\text{dim}\, V-1}{d} = 
\frac{P_{n-1}(d)+\text{dim}\, W-1}{d}\  . 
\]
This inequality is equivalent to $\text{dim}\, W - 1 < (d-1)P_{n-1}(d)$ which 
is equivalent to $\text{dim}\, W \leq (d-1)P_{n-1}(d)$. But $\text{dim}\, W 
\leq P_n(d-1)$ and it is easy to check that $P_n(d-1) \leq 
(d-1)P_{n-1}(d)$. 

(iii) If $u=X_n$ then $(V:u)=W$ and we have already shown that: 
\[
\frac{\text{dim}\, W-1}{d-1} < \frac{\text{dim}\, V-1}{d}\  . 
\]
\end{proof} 
 
In the next three remarks we show that, using the recursive constructions 
from the above proof of Theorem 4, one can make this proof independent of 
the constructions of Costa, Macias Marques and Mir\'{o}-Roig \cite{cmm}. 

\vskip3mm 
{\bf 3.1. Remark.}\quad 
\textit{Assume that} $n=1$ \textit{and consider two integers} $d \geq 1$ 
\textit{and} $2 \leq m \leq d+1 = P_1(d)$. \textit{Consider, also, an 
integer} $c \geq 1$ \textit{such that} $(m-1)c \leq d$. \textit{Then 
there exists an} $m$-\textit{dimensional b.p.f. monomial subspace} $V$ 
\textit{of} $S_d$ \textit{such that}, $\forall \, 1 \leq e \leq d-1$, 
$\forall \, u \in S_e$ \textit{monomial}: 
\[
\frac{\text{dim}(V:u)-1}{d-e} \leq \frac{1}{c}\  .
\]
 
\begin{proof} 
We take $V = kX_0^d + \sum_{p=0}^{m-2}kX_0^{pc}X_1^{d-pc}$. 

(i) If $u=X_0^iX_1^{e-i}$ with $i<e$ then: 
\[
\text{dim}(V:u) = \text{card}\{p\leq m-2\, \vert \, i\leq pc,\  e-i\leq 
d-pc\} = \text{card}\{p\leq m-2\, \vert \, i\leq pc \leq d-e+i\}
\]
hence $\text{dim}(V:u)-1 \leq (d-e)/c$. 

(ii) If $u=X_0^e$ then $\text{dim}(V:u) = 1+\text{card}\{p\leq m-2\, 
\vert \, e\leq pc\}$ hence: 
\[
\text{dim}(V:u)-1 = m-1-\lceil e/c \rceil \leq 
d/c-e/c = (d-e)/c\  . 
\]
\end{proof}         

{\bf 3.2. Remark.}\quad 
\textit{Assume that} $n=2$ \textit{and consider two integers} $d \geq 1$ 
\textit{and} $3 \leq m \leq 2d+1$. \textit{Then there exists an} 
$m$-\textit{dimensional b.p.f. monomial subspace} $V$ \textit{of} $S_d$ 
\textit{which satisfies the non-strict inequalities from Corollary 2.6}. 
\vskip3mm 

\begin{proof} 
We consider four integers (to be specified later) $m_1 \geq 2$, $m_2 \geq 2$, 
$c_1 \geq 1$ and $c_2 \geq 1$ which are supposed to satisfy the relations: 
\begin{gather}
m-1 = (m_1-1) + (m_2-1) \tag{1}\\ 
\frac{d}{m-1} \leq c_1 \leq \frac{d}{m_1-1},\  
\frac{d}{m-1} \leq c_2 \leq \frac{d}{m_2-1}\  . \tag{2}
\end{gather} 
Then we consider the monomial subspaces: 
\begin{gather*} 
V_{01} := kX_0^d + {\textstyle \sum_{p=0}^{m_1-2}}kX_0^{pc_1}X_1^{d-pc_1},\\ 
V_{02} := kX_0^d + {\textstyle \sum_{p=0}^{m_2-2}}kX_0^{pc_2}X_2^{d-pc_2},\\
V := V_{01} + V_{02}.
\end{gather*} 
Since $V_{01} \cap V_{02} = kX_0^d$, $V$ is $m$-dimensional. We assert that 
$V$ satisfies the non-strict inequalities from Corollary 2.6. Indeed, let 
$1 \leq e \leq d-1$ be an integer and $u \in S_e$ a monomial. 

(i) If $X_1X_2$ divides $u$ then $(V:u) = (0)$. 

(ii) If $u \in k[X_0,X_i]_e \setminus \{X_0^e\}$, $i = 1, 2$, 
then $(V:u) = (V_{0i}:u)$ hence, from Remark 3.1.: 
\[
\frac{\text{dim}(V:u)-1}{d-e} \leq \frac{1}{c_i} \leq \frac{m-1}{d} = 
\frac{\text{dim}\, V-1}{d}\  .
\]

(iii) If $u = X_0^e$ then, using the computations from the last part  
of the proof of Remark 3.1., one has: 
\begin{gather*} 
\text{dim}(V:u)-1 = (\text{dim}(V_{01}:u)-1) + (\text{dim}(V_{02}:u)-1) =\\ 
= (m_1-1-\lceil e/c_1\rceil ) + (m_2-1-\lceil e/c_2\rceil ) = 
m-1-\lceil e/c_1\rceil - \lceil e/c_2\rceil . 
\end{gather*}
The inequality we would like to prove: 
\[
\frac{\text{dim}(V:u)-1}{d-e} \leq \frac{\text{dim}\, V-1}{d} = 
\frac{m-1}{d} 
\]
is, consequently, equivalent to: 
\[
\left\lceil \frac{e}{c_1}\right\rceil + \left\lceil \frac{e}{c_2}\right\rceil 
\geq e\cdot \frac{m-1}{d}\  .
\]
But the last inequality is fulfilled because: 
\[
\frac{1}{c_1} + \frac{1}{c_2} \geq \frac{m_1-1}{d} + \frac{m_2-1}{d} = 
\frac{m-1}{d}\  .
\]
  
We notice, for further use, that if, moreover, one of the following two 
sets of stronger conditions: 
\begin{gather} 
\frac{d}{m-1} < c_1 < \frac{d}{m_1-1},\  
\frac{d}{m-1} < c_2 \leq \frac{d}{m_2-1} \tag{3} \\
\frac{d}{m-1} < c_1 \leq \frac{d}{m_1-1},\  
\frac{d}{m-1} < c_2 < \frac{d}{m_2-1}\  . \tag{4}
\end{gather} 
is fulfilled then the subspace $V$ constructed above satisfies the 
{\it strict} inequalities from Corollary 2.6. 

Finally, it remains to specify the integers $m_1$, $m_2$, $c_1$ and $c_2$.   

{\bf Case 1:} $3 \leq m \leq d+1$ \textit{and} $m$ \textit{odd}. 

Assume that $m = 2t+1$ for some $t \geq 1$. We take $m_1 = m_2 = t+1$ 
and $c_1 = c_2 = \lfloor d/t\rfloor$. Since $d \geq m-1 = 2t$ one has  
$d/(m-1) = d/2t \leq d/t - 1 < 
\lfloor d/t\rfloor $  
hence the relations (1) and (2) are fulfilled. 

{\bf Case 2:} $3 \leq m \leq d+1$ \textit{and} $m$ \textit{even}. 

Assume that $m = 2t+2$ for some $t \geq 1$. We take $m_1 = t+2$, 
$m_2 = t+1$, $c_1 = \lfloor d/(t+1)\rfloor$ and $c_2 = \lfloor d/t\rfloor$. 
Since $d \geq m-1 = 2t+1$ one has: 
\[
\frac{d}{m-1} = \frac{d}{2t+1} \leq \frac{d}{t+1} - \frac{t}{t+1} \leq 
\left\lfloor \frac{d}{t+1}\right\rfloor 
\]
hence the relations (1) and (2) are fulfilled. 

{\bf Case 3:} $d+2 \leq m \leq 2d+1$. 

We take $m_1 = d+1$, $m_2 = m-d$ and $c_1 = c_2 = 1$. The relations (1) and 
(2) are obviously fulfilled. In this case, the subspace $V$ constructed 
above equals $S^{\prime}_d + k[X_0,X_2]_l\cdot X_2^{d-l}$, where  
$l = m-d-2$, and appears in the proof of \cite{cmm}, Proposition 3.3.   
\end{proof}
\vskip3mm 

{\bf 3.3. Remark.}\quad 
\textit{Assume that} $n=2$. \textit{Then}, \textit{for any two integers} 
$d \geq 2$ \textit{and} $2d+2 \leq m \leq P_2(d)$ \textit{there exists an} 
$m$-\textit{dimensional b.p.f. monomial subspace} $V$ \textit{of} $S_d$ 
\textit{which satisfies the strict inequalities from Corollary 2.6}. 
\vskip3mm 

\begin{proof} 
We use induction on $d$. In the case $d=2$ one takes $V = S_d$. Let us treat 
now the {\it induction step} $(d-1) \rightarrow d \geq 3$. Since 
$P_2(d) - d - 1 = P_2(d-1)$, we can write $m = d + 1 + l$, with 
$d + 1 \leq l \leq P_2(d-1)$. From Remark 3.2. and from the induction 
hypothesis, there exists an $l$-dimensional b.p.f. monomial subspace 
$W$ of $S_{d-1}$ such that, $\forall \, 1 \leq e \leq d-2$, 
$\forall \, u \in S_e$ monomial: 
\[
\frac{\text{dim}(W:u)-1}{d-1-e} \leq \frac{\text{dim}\, W-1}{d-1}\  .\tag{1} 
\]
 We take $V := S^{\prime}_d + WX_2$. Consider an integer $1 \leq e \leq d-1$ 
and a monomial $u \in S_e$. 

(i) If $u \in S^{\prime}_e$ then $(V:u) = S^{\prime}_{d-e} + (W:u)X_2$. 
Since $\text{dim}\, S^{\prime}_{d-e} = d-e+1$ and $\text{dim}\, S^{\prime}_d 
= d+1$, the strict inequality from Corollary 2.6. is equivalent, in this 
case, to: 
\[
\frac{\text{dim}(W:u)}{d-e} < \frac{\text{dim}\, W}{d}\  . \tag{2} 
\]

Now, when $e=d-1$ one has $\text{dim}(W:u) \leq 1$, $d-e = 1$ and 
$\text{dim}\, W \geq d+1$ hence (2) is fulfilled. When $1 \leq e \leq d-2$, 
it follows from (1) that:
\[
\text{dim}(W:u) \leq \frac{d-e-1}{d-1}\text{dim}\, W + \frac{e}{d-1}\  .
\]
But, since $\text{dim}\, W \geq d+1$, one checks easily that: 
\[
\frac{d-e-1}{d-1}\text{dim}\, W + \frac{e}{d-1} < \frac{d-e}{d}
\text{dim}\, W 
\]
hence (2) is fulfilled. 

(ii) If $u \in S_{e-1}X_2$ and $e \geq 2$ one uses the argument from the 
Case 3(ii) of the proof of Theorem 4. 

(iii) If $u=X_2$ one uses the argument from the Case 3(iii) of the proof of 
Theorem 4. 
\end{proof}
\vskip3mm 

{\bf 3.4. Remark.}\quad 
\textit{If} $n=2$  \textit{and} $3 \leq m \leq 2d+1$ \textit{then}, 
\textit{except for the cases} $m = d+1$ \textit{and} $m = 2d+1$,   
\textit{the construction from the proof of Remark 3.2. can be used to 
produce} $m$-\textit{dimensional b.p.f. monomial 
subspaces} $V$ \textit{of} $S_d$ \textit{which satisfy the strict 
inequalities from Corollary 2.6}. 
\vskip3mm 

\begin{proof}    
We use the notation from the proof of Remark 3.2. 

{\bf Case 1:} $3 \leq m \leq d$ \textit{and} $m$ \textit{odd}. 

Assume that $m = 2t+1$ for some $t \geq 1$. If $t$ does not divide $d$ then 
the integers $m_1$, $m_2$, $c_1$, $c_2$ defined in the Case 1 of the proof of 
Remark 3.2. satisfy the relations (1), (3) and (4) 
from the proof of that remark hence 
the monomial subspace $V$ of $S_d$ considered at the beginning of that proof  
satisfies the strict inequalities from Corollary 2.6. 

If $t$ divides $d$ we take $m_1 = m_2 = t+1$ and $c_1 = c_2 = d/t-1$. Since 
$d \geq m = 2t+1$ it follows that 
$d/(m-1) = d/2t < d/t-1$ hence $m_1$, $m_2$, $c_1$, $c_2$  
satisfy the relations (1), (3) and (4) from the proof of Remark 3.2. 

{\bf Case 2:} $3 \leq m \leq d$ \textit{and} $m$ \textit{even}. 

Assume that $m = 2t+2$ for some $t \geq 1$. Since $d \geq m = 2t+2$ one has: 
\[
\frac{d}{m-1} = \frac{d}{2t+1} < \frac{d}{t+1} - \frac{t}{t+1} \leq 
\left\lfloor \frac{d}{t+1}\right\rfloor \  .
\] 

If $t$ does 
not divide $d$ then the integers $m_i$, $c_i$ defined in the Case 2 of 
the proof of Remark 3.2 satisfy the relations (1) and (4) from that proof. 

Now, if $t$ divides $d$ we take $m_1 = t+2$, $m_2 = t+1$, 
$c_1 = \lfloor d/(t+1)\rfloor $ and $c_2 = d/t-1$. 
Since $d \geq m = 2t+2$  we have  
$d/(m-1) = d/(2t+1) < d/2t < d/t-1$. The integers $m_i$, $c_i$ satisfy in 
this case the relations (1) and (4) from the proof of Remark 3.2. 

{\bf Case 3:} $d+2 \leq m\leq 2d$. 

The integers $m_i$, $c_i$ defined in the Case 3 of the proof of Remark 3.2. 
satisfy the relations (1) and (4) from that proof.  
\end{proof} 
\vskip3mm 
 
In the cases $n=2$, $m=d+1$ and $n=2$, $m=2d+1$ one needs different 
constructions if one wants to get monomial subspaces satisfying the strict 
inequalities from Corollary 2.6. 
For example, if $n = 2$, $m = d+1$ and 
$d \geq 6$ one may take $V = V^{\prime} + kX_1^aX_2^b$, where $V^{\prime}$ is 
the $d$-dimensional b.p.f. monomial subspace of $S_d$ constructed in the 
proof of Remark 3.2. and $a \geq 3$, $b \geq 3$ are two integers with 
$a + b = d$. 

Indeed, let $1\leq e \leq d-1$ be an integer and $u \in S_e$ a monomial. 
If $u$ does not divide $X_1^aX_2^b$ then $(V:u) = (V^{\prime}:u)$ hence: 
\[
\frac{\text{dim}(V:u)-1}{d-e} = \frac{\text{dim}(V^{\prime}:u)-1}{d-e} 
\leq \frac{\text{dim}\, V^{\prime}-1}{d} < \frac{\text{dim}\, V-1}{d}\  . 
\]
If $X_1X_2$ divides $u$ then $\text{dim}(V:u)-1 \leq 0$. It remains to 
consider the cases $u = X_1^e$, $e \leq a$ and $u = X_2^e$, $e \leq b$. In 
the first case, noticing that the integers $c_1$, $c_2$ used in the 
construction of $V^{\prime}$ are $\geq 2$, one has: 
\begin{gather*} 
\frac{\text{dim}(V:u)-1}{d-e} = 
\frac{\text{dim}(V^{\prime}_{01}:u) + 1 -1}{d-e}  
\leq \frac{1}{2} + \frac{1}{d-e}\leq \\ 
\leq \frac{1}{2} + \frac{1}{d-a} = \frac{1}{2} + \frac{1}{b} < 1 = 
\frac{\text{dim}\, V-1}{d}\  .
\end{gather*}
For $u = X_2^b$ one uses a similar argument. 

If $n = 2$, $m = d+1$ and $2 \leq d \leq 5$ one may take, respectively: 
\begin{gather*} 
V = kX_0^2 + kX_1^2 + kX_2^2,\  V = V^{\prime} + kX_0X_1X_2,\\ 
V = V^{\prime} + kX_0X_1X_2^2,\  V = V^{\prime} + kX_0X_1^2X_2^2,  
\end{gather*} 
where $V^{\prime}$ is as above.   

Analogously, if $n = 2$, $m = 2d+1$ and $d \geq 4$ one may take 
$V = V^{\prime} + kX_1^aX_2^b$, where $V^{\prime}$ is the $2d$-dimensional  
b.p.f. monomial subspace of $S_d$ constructed in the proof of Remark 3.2. and 
$a \geq 2$, $b \geq 2$ are two integers with $a + b = d$. 

Finally, if $n = 2$, $d = 3$ and $m = 2d+1 = 7$ one may take: 
\[
V = kX_0^3 + kX_0X_1^2 + kX_1^3 + kX_0^2X_2 + kX_0X_1X_2 + kX_1X_2^2 + 
kX_2^3 . 
\]

\end{document}